\documentclass[reqno,a4paper,11pt]{amsart}
\usepackage{hyperref}
\usepackage{amsthm,amssymb}
\usepackage{eepic}
\usepackage{enumerate}
\usepackage{xspace}
\theoremstyle{definition}
\newtheorem{definition}{Definition}
\newtheorem{example}[definition]{Example}
\theoremstyle{remark}

\theoremstyle{plain}
\newtheorem{lemma}[definition]{Lemma}

\newtheorem{theorem}[definition]{Theorem}
\newtheorem{cor}[definition]{Corollary}
\newtheorem{corollary}[definition]{Corollary}

\newcommand{\set}[1]{\left\{{#1}<\right\}}
\newcommand{\vek}[1]{\mathbf{#1}}

\newcommand\setsuchas[2]{\left\{\,{#1}\,\vrule\,{#2}\,\right\}}

\newcommand{\Nat}{{\mathbb{N}}}
\newcommand{\Z}{{\mathbb{Z}}}
\newcommand{\C}{{\mathbb{C}}}

\newcommand{\Q}{{\mathbb{Q}}}

\newcommand{\PartitionOf}[2]{#1 \vdash #2}
\newcommand{\NumParts}[1]{\langle #1 \rangle}
\newcommand{\w}[1]{{\lvert #1 \rvert}}
\newcommand{\card}[1]{{\lvert #1 \rvert}}
\newcommand{\Partitions}{\mathcal{P}}

\newcommand{\pseudoconcave}{super-concave\xspace}
\newcommand{\concave}{concave\xspace}

\begin{document}
\title{Enumeration of concave integer partitions}
\author{Jan Snellman}
\address{Department of Mathematics\\
Stockholm University\\
SE-10691 Stockholm, Sweden}
\email{Jan.Snellman@math.su.se}
\author{Michael Paulsen}
\begin{abstract}
  An integer partition \(\PartitionOf{\lambda}{n}\) corresponds, via
  its Ferrers diagram, to an artinian monomial ideal \(I \subset
  \C[x,y]\) with \(\dim_\C \C[x,y]/I = n\). If \(\lambda\) corresponds
  to an integrally closed ideal we call it \emph{\concave}. We study
  generating functions for the number of \concave partitions,
  unrestricted or with at most \(r\) parts.
\end{abstract}

\maketitle

\begin{section}{\concave partitions}
  By an \emph{integer partition} \(\lambda = (\lambda_1, \lambda_2,
  \lambda_3,\dots)\)
  we mean a weakly decreasing sequence of non-negative integers, all
  but finitely many of which are zero. The non-zero elements are
  called the \emph{parts} of the partition. When writing a partition,
  we often will only write the parts; thus \((2,1,1,0,0,0,\dots)\) may
  be written as \((2,1,1)\).

  We write
  \(r=\NumParts{\lambda}\)
  for the number of parts of \(\lambda\), and \(n=\w{\lambda}=\sum_{i}
  \lambda_i\); equivalently, we write \(\PartitionOf{\lambda}{n}\)
  if \(n=\w{\lambda}\). The set of all partitions is denoted by
  \(\Partitions\), and the set of partitions of \(n\) by
  \(\Partitions(n)\). We put
  \(\card{\Partitions(n)}=p(n)\). By subscripting any of the above
  with \(r\) we
  restrict to partitions with at most \(r\) parts.

  We will use the fact that \(\Partitions\) forms a monoid
  under component-wise addition.

  For an integer partition \(\PartitionOf{\lambda}{n}\) we define its
  \emph{Ferrers diagram} \(F(\lambda) = \setsuchas{(i,j) \in \Nat^2}{i
    < \lambda_{j+1}}\). In figure~\ref{fig:exejkonk} the black dots
  comprise the Ferrers diagram of the partition \(\mu=(4,4,2,2)\).

Then \(F(\lambda)\) is a finite \emph{order ideal}
  in the partially ordered set \((\Nat^2,\le)\), where \((a,b) \le
  (c,d)\) iff \(a \le c\) and \(b \le d\). In fact, integer partitions
  correspond precisely to finite order ideals in this poset.

  The complement \(I(\lambda)=\Nat^2 \setminus F(\lambda)\) is a
  monoid ideal in the additive monoid \(\Nat^2\). Recall that for a
  monoid ideal \(I\) the \emph{integral closure} \(\bar{I}\) is
  \begin{equation}
    \label{eq:intcl}
    \setsuchas{\vek{a}}{ \ell \vek{a} \in I \text{ for some } \ell \in
      \Z_+}
  \end{equation}
  and that \(I\)
  is \emph{integrally closed} iff it is equal to its integral closure.

  \begin{definition}\label{def:ic}
    The integer partition \(\lambda\) is \emph{\concave} iff
    \(I(\lambda)\) is integrally closed. We denote by
    \(\bar{\lambda}\) the unique partition such that
    \(I(\bar{\lambda})= \overline{I(\lambda)}\).
  \end{definition}

  Now let \(R\) be the complex monoid ring of \(\Nat^2\).
  We identify \(\Nat^2\) with the set of commutative monomials in the
  variables \(x,y\), so that
  \(R\simeq \C[x,y]\). Then a monoid ideal \(I \subset \Nat^2\)
  corresponds to the monomial ideal \(J\) in \(R\) generated by the monomials
  \(\setsuchas{x^iy^j}{(i,j) \in I}\). Furthermore, since the monoid
  ideals of the form \(I(\lambda)\) are precisely those with finite
  complement to \(\Nat^2\), those monoid ideals will correspond to
  monomial ideals \(J \subset R\) such that \(R/J\) has a finite
  \(\C\)-vector space basis (consisting of images of those monomials
  not in \(J\)). By abuse of notation, such monomial ideals are called
  \emph{artinian}, and the \(\C\)-vector space dimension of \(R/J\) is
  called the \emph{colength} of \(J\).

  We get in this way a bijection between
  \begin{enumerate}
  \item integer partitions of \(n\),
  \item order ideals in \((\Nat^2,\le)\) of cardinality \(n\),
  \item monoid ideals in \(\Nat^2\) whose complement has cardinality
    \(n\), and
  \item monomial ideals in \(R\) of colength \(n\).
 \end{enumerate}

 Recall that if \(\mathfrak{a}\) is an ideal in the
 commutative unitary ring \(S\), then the \emph{integral closure}
 \(\bar{\mathfrak{a}}\) consists of all \(u \in S\) that fulfill some
 equation of the form
 \begin{equation}
   \label{eq:iceq}
   s^n + b_1s^{n-1} + \dots + b_0, \qquad b_i \in \mathfrak{a}^i
 \end{equation}
 Then \(\mathfrak{a}\) is always contained in its integral closure,
 which is an ideal. The ideal \(\mathfrak{a}\) is said to be
 \emph{integrally closed} if it coincides with its integral closure.

 For the special case \(S=R\), we have that the integral closure of a
 monomial ideal is again a monomial ideal, and that the latter monomial
 ideal corresponds to the integral closure of the monoid ideal
 corresponding to the former monomial ideal. Hence, we have a
 bijection between

 \begin{enumerate}
 \item \concave integer partitions of \(n\),
 \item integrally closed monoid ideals in \(\Nat^2\) whose complements
   have cardinality
   \(n\), and
 \item integrally closed monomial ideals in \(R\) of colength \(n\).
 \end{enumerate}

 Fr{\"o}berg and Barucci \cite{Froeberg:Finco}
 studied the growth of the number of ideals of colength \(n\)
 in certain rings, among them local noetherian rings of dimension 1.
 Studying the growth of the number of monomial
 ideals of colength \(n\) in \(R\) is, by the above, the same as
 studying the partition function \(p(n)\). In this article, we will
 instead study the growth of the number of integrally closed monomial
 ideals in \(R\), that is, the number of \concave partitions of \(n\).

\begin{figure}[bt]
\setlength{\unitlength}{0.75 cm}
\begin{picture}(15,5)
  \put(0,0){\line(1,0){5.5}}
  \put(0,0){\line(0,1){5}}
 
  \multiput(-0.12,-0.135)(1,0){4}{$\bullet$}
  \multiput(-0.12,0.855)(1,0){4}{$\bullet$}
  \multiput(-0.12,1.855)(1,0){2}{$\bullet$}
  \multiput(-0.12,2.855)(1,0){2}{$\bullet$}
  \multiput(-0.196,3.87)(1,0){6}{$\times$}
  \multiput(1.814,2.87)(1,0){4}{$\times$}
  \multiput(1.814,1.87)(1,0){4}{$\times$}
  \multiput(3.814,0.87)(1,0){2}{$\times$}
  \multiput(3.814,-0.13)(1,0){2}{$\times$}

  \put(10,0){\line(1,0){5.5}}
  \put(10,0){\line(0,1){5}}
 
  \multiput(9.88,-0.135)(1,0){4}{$\bullet$}
  \multiput(9.88,0.855)(1,0){3}{$\bullet$}
  \multiput(9.88,1.855)(1,0){2}{$\bullet$}
  \multiput(9.88,2.855)(1,0){1}{$\bullet$}
  \multiput(9.804,3.87)(1,0){6}{$\times$}
  \multiput(10.814,2.87)(1,0){5}{$\times$}
  \multiput(11.814,1.87)(1,0){4}{$\times$}
  \multiput(12.814,0.87)(1,0){3}{$\times$}
  \multiput(13.814,-0.13)(1,0){2}{$\times$}

\end{picture}
\caption{\(\mu\) and \(\bar{\mu}\)} \label{fig:exejkonk}
\end{figure}
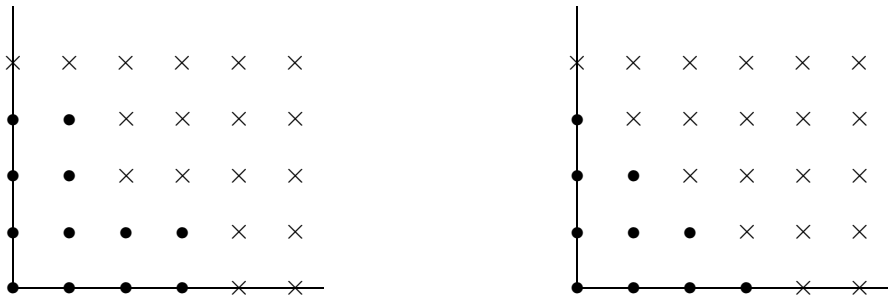 







\end{section}

\begin{section}{Inequalities defining \concave partitions}
  It is in general a hard problem to compute the integral closure of
  an ideal in a commutative ring. However, for monomial ideals in a
  polynomial ring, the following theorem, which can be found in
  e.g. \cite{Ebud:View}, makes the problem feasible.
  \begin{theorem}\label{thm:convex}
    Let \(I \subset \Nat^2\) be a monoid ideal, and regard \(\Nat^2\)
    as a subset of \(\Q^2\) in the natural way. Let
    \(\mathrm{conv}_\Q(I)\) denote the convex hull of \(I\) inside
    \(\Q^2\).
    Then the integral closure of \(I\) is given by
    \begin{equation}
      \label{eq:latticep}
      \mathrm{conv}_\Q(I) \cap \Nat^2
    \end{equation}
  \end{theorem}

\begin{example}
The partition \(\mu=(4,4,2,2)\) corresponds to the monoid ideal
\(((0,4),(2,2),(4,0))\), which has  integral closure
\(((0,4),(1,3),(2,2),(3,1),(4,0))\).
It follows that \(\overline{\mu}=(4,3,2,1)\).
In figure~\ref{fig:exejkonk} we have drawn the lattice points belonging
to \(F(\mu)\) as dots, and the lattice points belonging to
\(I(\lambda)\) as crosses.
\end{example}

The above theorem gives the following characterization of
\concave partitions:
\begin{lemma}\label{lemma:concaveChar}
  Let \(\lambda=(\lambda_1,\lambda_2,\lambda_3,\dots)\) be a
  partition. Then \(\lambda\) is \concave iff
  for all positive integers \(i < j < k\),
  \begin{equation}
    \label{eq:konkav}
    \lambda_j <  1+ \lambda_i \frac{k-j}{k-i} + \lambda_k
    \frac{j-i}{k-i}
  \end{equation}
  or, equivalently, if
  \begin{equation}
    \label{eq:konkavhom}
    \lambda_i (j-k)  + \lambda_j (k-i) + \lambda_k (i-j) < k-i
  \end{equation}
\end{lemma}

\end{section}

\begin{section}{Generating functions for \pseudoconcave partitions}
We will enumerate \concave partitions by considering another class of
partitions which is more amenable to enumeration, yet is close to
that of \concave partitions.
\begin{definition}\label{def:almostconcave}
  Let \(\lambda=(\lambda_1,\lambda_2,\lambda_3,\dots)\) be a
  partition. Then \(\lambda\) is \emph{\pseudoconcave} iff
  for all positive integers \(i < j < k\),
  \begin{equation}
    \label{eq:pkonkav}
    \lambda_i (j-k)  + \lambda_j (k-i) + \lambda_k (j-i) \le 0
  \end{equation}
\end{definition}

The reader should note that it is actually a
\emph{stronger} property to be \pseudoconcave than to be \concave. 
Unlike the latter property, it is not
necessarily preserved by conjugation: the partition \((2)\) is
\pseudoconcave, hence \concave, but its conjugate \((1,1)\) is \concave
but not \pseudoconcave.

\begin{theorem}\label{thm:pseudoconcave}
  Let \(\lambda=(\lambda_1,\lambda_2,\lambda_3,\dots)\) be a
  partition, and let \(\mu =(\mu_1,\mu_2,\mu_3,\dots)\) be its
  conjugate, so that \(\card{\setsuchas{j}{\mu_j=i}} = \lambda_i -
  \lambda_{i+1}\) for all \(i\).
  Then the following are equivalent:
  \begin{enumerate}[(i)]
  \item \label{it:ps}
    \(\lambda\) is \pseudoconcave,
  \item \label{it:sv}
    for all positive \(\ell\),
    \begin{equation}
      \label{eq:asc}
            -\lambda_\ell +2\lambda_{\ell+1} - \lambda_{\ell+2} \le 0
    \end{equation}
  \item \label{it:asc}
    for all positive \(\ell\),
    \begin{equation}
      \label{eq:asc2}
            \lambda_{\ell+1} -\lambda_\ell \ge \lambda_{\ell+2} -
            \lambda_{\ell+1}
    \end{equation}
  \item  \label{it:sfp}
    \(\card{\setsuchas{k}{\mu_{k}=i}} \geq
      \card{\setsuchas{k}{\mu_{k}=j}}\) whenever \(i \leq j\).
  \end{enumerate}
\end{theorem}
\begin{proof} \eqref{it:ps} \(\iff\) \eqref{it:sv}:
  Let \(\vek{e}_i\) be the vector with 1 in the \(i\)'th coordinate
  and zeros elsewhere, let \(\vek{f}_j=-\vek{e}_j + 2\vek{e}_{j+1}
  -\vek{e}_{j+2}\), and let \(\vek{t}_{i,j,k} =
    (j-k) \vek{e}_i  +  (k-i)\vek{e}_j  + (j-i)\vek{e}_k\).
  Clearly, \eqref{eq:pkonkav} is equivalent with
  \(\vek{t}_{i,j,k} \cdot \lambda \le 0\), and
  \eqref{eq:asc} is equivalent with
  \(\vek{f}_{j} \cdot \lambda \le 0\).
  We have that \(\vek{f}_{\ell} = \vek{t}_{\ell,\ell+1,\ell+2}\).
  Conversely, we claim that \(\vek{t}_{i,j,k}\)
   is a positive linear combination of different \(\vek{f}_{\ell}\).
  From this claim, it follows
  that if \(\lambda\) fulfills \eqref{eq:asc} for all \(\ell\)
  then \(\lambda\) is \pseudoconcave.

  We can without loss of generality assume that \(i=1\).
  Then it is easy to verify that
  \begin{equation}
    \label{eq:poslincomb}
    \vek{t}_{1,j,k} = \sum_{\ell=1}^{j-2} \ell (k-j)\vek{f}_{\ell}
    + \sum_{\ell=j-1}^{k-2} \ell (j-1)(k-\ell - 1) \vek{f}_{\ell}
  \end{equation}

  \eqref{it:sv} \(\iff\) \eqref{it:asc}    \(\iff\)  \eqref{it:sfp} :
  This is obvious.
\end{proof}

The \emph{difference operator} \(\Delta\) is defined on partitions by
\begin{equation}
  \label{eq:delta}
  \Delta(\lambda_1,\lambda_2,\lambda_3,\dots) = (\lambda_1 -
  \lambda_2, \lambda_2 - \lambda_3, \lambda_3-\lambda_4,\dots)
\end{equation}
We get that the \emph{second order difference operator} \(\Delta^2\)
is given by
\begin{multline}
  \label{eq:delta2}
    \Delta^2(\lambda_1,\lambda_2,\lambda_3,\dots) =   
    \Delta(\Delta(\lambda_1,\lambda_2,\lambda_3,\dots)) = \\
    = (\lambda_1 - 2\lambda_2 + \lambda_3, \lambda_2 - 2\lambda_3 +
    \lambda_4, \lambda_3 - 2\lambda_4 + \lambda_5, \dots)
\end{multline}
\begin{corollary}
  The \pseudoconcave partitions are precisely those with non-negative
  second differences.
\end{corollary}

\begin{definition}
    Let \(p_{sc}(n)\) denote the number of \pseudoconcave partitions of
  \(n\), and \(p_{sc}(n,r)\) denote the number of \pseudoconcave
  partitions of
  \(n\) with at most \(r\) parts.
  Let similarly \(p_{c}(n)\) and \(p_{c}(n,r)\)
denote the number of \pseudoconcave partitions of \(n\),
and the number of \pseudoconcave partitions of
\(n\) with at most \(r\) parts, respectively.
  For a partition \(\lambda=(\lambda_1,\lambda_2,\dots)\) let
  \(\vek{x}^\lambda= x_1^{\lambda_1}  x_2^{\lambda_2} \cdots\),
  and define
  \begin{equation}
    \label{eq:PS}
    \begin{split}
    PS(\vek{x}) &= \sum_{\lambda \text{ \pseudoconcave}} \vek{x}^\lambda\\
    PS_r(x_1,\dots,x_r) &= PS(x_1,x_2,\dots,x_r,0,0,0,\dots)
    = \sum_{\substack{\lambda \text{ \pseudoconcave}\\
        \lambda_{r+1}=0}} \vek{x}^\lambda \\
    PC(\vek{x}) &= \sum_{\lambda \text{ \concave}} \vek{x}^\lambda\\
    PC_r(x_1,\dots,x_r) &= PC(x_1,x_2,\dots,x_r,0,0,0,\dots)
    = \sum_{\substack{\lambda \text{ \concave}\\
        \lambda_{r+1}=0}} \vek{x}^\lambda
    \end{split}
  \end{equation}
\end{definition}

Partitions with non-negative second differences have been studied by
Andrews \cite{Andrews:MM2}, who proved that there are as many such
partitions of \(n\) as 
there are partitions of \(n\) into triangular numbers.

Canfield et al \cite{CCH:RandomPart} have studied partitions with
non-negative \(m\)'th differences. Specialising their results to the
case \(m=2\), we conclude:

\begin{theorem}\label{thm:nonnegsecdiff}
  Let \(n,r\) be  denote positive integers.

  \begin{enumerate}[(i)]
  \item There is a bijection between
    partitions of \(n\) into triangular numbers and \pseudoconcave
    partitions.
  \item     The multi-generating function for \pseudoconcave partitions is
    given by
    \begin{equation}
      \label{eq:psgen}
      \begin{split}
      PS(\vek{x}) &=  \frac{1}{ \prod_{i=1}^\infty \left(1 - \prod_{j=1}^i
        x_j^{1+i-j}\right)}\\
      &=
      1  +
      x_{{1}} +
      {x_{{1}}}^{2} +
      {x_{{1}}}^{3} +
      {x_{{1}}}^{4} +
      {x_{{1}}}^{2}x_{{2}} +
      {x_{{1}}}^{5} +
      {x_{{1}}}^{4}x_{{2}} +
      {x_{{1}}}^{3}x_{{2}}     + \dots
      \end{split}
      \end{equation}
    \item     The multi-generating function for \pseudoconcave partitions with at
    most \(r\) parts is
    given by
    \begin{equation}
      \label{eq:psgenrp}
      PS_r(x_1,x_2,\dots,x_r) =
\frac{1}{ \prod_{i=1}^r \left(1 - \prod_{j=1}^i
        x_j^{1+i-j}\right)}
    \end{equation}
  \item
    The generating function for \pseudoconcave partitions is
    \begin{equation}
      \label{eq:psgenspec}
      PS(t) = \sum_{n=0}^\infty p_{sc}(n)t^n =
      \prod_{i=1}^\infty \frac{1}{ 1 - t^{\frac{i(i+1)}{2}}}
    \end{equation}
    and the one for \pseudoconcave partitions with at most \(r\) parts
    is
    \begin{equation}
      \label{eq:psgenspecrp}
      PS_r(t) = \sum_{n=0}^\infty p_{sc}(n,r)t^n =
      \prod_{i=1}^r \frac{1}{ 1 - t^{\frac{i(i+1)}{2}}}
    \end{equation}

    \item       The proportion of \pseudoconcave partitions with
      at most \(r\) parts among all partitions with at most \(r\) parts
      is
    \begin{equation}
      \label{eq:prop}
      \frac{r!}{\prod_{i=1}^r \frac{i(i+1)}{2}}.
    \end{equation}

  \item As \(n \to \infty\),
      \begin{equation}
        \label{eq:psas}
        \begin{split}
         p_{sc}(n)  \sim c n^{-3/2} \exp(3C n^{1/3})\\
         C = 2^{-1/3} \left[\zeta (3/2)\Gamma(3/2)\right]^{2/3}, \quad
         c = \frac{\sqrt{3}}{12}\,\left ({\frac {C}{\pi }}\right )^{3/2}
        \end{split}
      \end{equation}
  \end{enumerate}
\end{theorem}
  The sequence \(\bigl(p_{sc}(n)\bigr)_{n=0}^\infty\) is identical to
  sequence 
\htmladdnormallink{A007294}{http://www.research.att.com/projects/OEIS?Anum=007294}
in OEIS \cite{Sloane}.
    We have submitted the sequences
    \(\bigl(p_{sc}(n,r)\bigr)_{n=0}^\infty\), for \(r=3,4\), in
    OEIS \cite{Sloane}, as 
\htmladdnormallink{A086159}{http://www.research.att.com/projects/OEIS?Anum=086159}
and 
\htmladdnormallink{A086160}{http://www.research.att.com/projects/OEIS?Anum=086160}.
The sequence for \(r=2\) was already in the database, as 
\htmladdnormallink{A008620}{http://www.research.att.com/projects/OEIS?Anum=008620}.

\begin{subsection}{Other apperances of \pseudoconcave partitions in
  the literature} 
  The bijection between partitions into triangular numbers and
  partitions with non-negative second difference  is mentioned in
  \htmladdnormallink{A007294}{http://www.research.att.com/projects/OEIS?Anum=007294}  
in OEIS \cite{Sloane}, together with a reference to Andrews
\cite{Andrews:MM2}.  That sequence has been contributed by  Mira
  Bernstein and Roland Bacher; 
  we thank Philippe Flajolet for drawing our attention to it.

Gert Almkvist \cite{Almkvist:partitions} gives an asymptotic
analysis of \(p_{sc}(n)\) which is finer than \eqref{eq:psas}.

Another derivation of the generating functions above can found in a
forthcoming paper ``Partition Bijections, a Survey'' \cite{Pak:Survey}
by Igor Pak. He observes that the set of 
\pseudoconcave  partitions with at most \(r\) parts consists of the
lattice points of the unimodular cone spanned by the vectors
\(v_0=(1,\dots,1)\) and \(v_i =(i-1,i-2,\dots,1,0,0,\dots)\) for \(1
\le i \le r\). 

Corteel and Savage \cite{PCInEq} calculate rational generating functions for
classes of partitions defined by linear homogeneous inequalities. This
applies to \pseudoconcave partitions, but not directly to concave partitions,
since the inequalities \eqref{eq:konkavhom} defining them are inhomogeneous.
\end{subsection}
\end{section}

    \begin{section}{Generating functions for \concave partitions}

    \begin{theorem}\label{thm:pcoco}
      Let \(r\) be a positive integer. Then
      \begin{equation}
        \label{eq:PCr}
        PC_r(x_1,\dots,x_r) =
        \frac{Q_r(x_1,\dots,x_r)}
{\prod_{i=1}^r \left( 1 - \prod_{j=1}^i
        x_j^{1+i-j}\right)}
      \end{equation}
      where \(Q_r(x_1,\dots,x_r)\) is a polynomial satisfying
      \begin{enumerate}[(i)]
      \item \(Q_r(x_1,\dots,x_r)\) has integer coefficients,
      \item \(Q_r(1,\dots,1)=1\),
      \item all exponent vectors of the monomials that occur
      in \(Q_r\) are weakly decreasing, and
      \item \(Q_r(x_1,\dots,x_r)=Q_{r+1}(x_1,\dots,x_r,0)\).
      \end{enumerate}
    Furthermore,
      \begin{equation}
        \label{eq:PC}
        PC(\vek{x}) =
        \frac{Q(\vek{x})}
        {\prod_{i=1}^\infty \left( 1 - \prod_{j=1}^i
            x_j^{1+i-j}\right)}
      \end{equation}
      where \(Q(\vek{x})\) is a formal power series with the property
      that for each \(\ell\),
      \(Q(x_1,\dots,x_\ell,0,0,\dots)=Q_\ell(x_1,\dots,x_\ell)\); in
      other words,
      \begin{displaymath}
        Q = 1 + \sum_{i=1}^\infty \left(Q_i - Q_{i-1}\right)
      \end{displaymath}
        \end{theorem}
    \begin{proof}
      Let \(A\) be the matrix with \(r\) columns whose rows consists
      of all  truncations of the vectors \(\vek{t}_{i,j,k}\)
      introduced in the proof of Theorem~\ref{thm:pseudoconcave},
      for \(i<j<k\), \(k<r+2\)
      For example, if \(r=3\) and if we order the 3-subsets of
      \(\set{1,2,3,4}\) lexicographically we get that
      \begin{displaymath}
        A =
        \begin{pmatrix}
          -1 & 2 & -1  \\
          -2 & 3 & 0  \\
          -1 & 0 & 3 \\
          0 & -1 & 2
        \end{pmatrix}
      \end{displaymath}
      Then a \pseudoconcave partition with at most \(r\) parts
      corresponds to a solution to
      \begin{equation}
        \label{eq:A}
        A\vek{z} \le \vek{0}, \qquad \vek{z} \ge \vek{0}
      \end{equation}
      whereas a \concave partition with at most \(r\) parts
      corresponds to a solution to
      \begin{equation}
        \label{eq:b}
        A\vek{z} \le \vek{b}, \qquad \vek{z} \ge \vek{0}
      \end{equation}
      where the entry of \(\vek{b}\) which corresponds to the row of
      \(A\) indexed by \((i,j,k)\) is \(i-k\).
      It follows from a theorem in Stanley's ``green book''
      \cite{Stanley:CombCom} that the multi-generating functions of
      these two solution sets have the same denominator, and that
      their numerator evaluates to the same value after substituting 1
      for each formal variable. 

      All monomials in
      \begin{displaymath}
        {\prod_{i=1}^r \left( 1 - \prod_{j=1}^i
        x_j^{1+i-j}\right)}
      \end{displaymath}
      have weakly decreasing exponent vectors, hence this is also true
      for \(PC_r(x_1,\dots,x_r)\).

      The assertion about \(PC(\vek{x})\) follows by passing to the
      limit.
    \end{proof}

    Our calculations indicate that
    \begin{equation}
      \label{eq:Qr}
      \begin{split}
      Q_1(\vek{x}) &= 1 \\
      Q_2(\vek{x}) & =1 + x_1x_2 - x_1^2x_2\\
      Q_3(\vek{x}) & =Q_2(\vek{x})  +  x_3\left({x_{{1}}}^{5}{x_{{2}}}^{3} -
        {x_{{1}}}^{4}{x_{{2}}}^{3} - 2\,{x_{{1}}}^{3}
{x_{{2}}}^{2} + {x_{{1}}}^{2}{x_{{2}}}^{2} + x_{{1}}x_{{2}}
 \right)\\
      \end{split}
    \end{equation}

    \begin{cor}\label{cor:convenum}
      \begin{enumerate}[(i)]
      \item 
      The generating function for \concave partitions with at most
      \(r\) parts is given by
      \begin{equation}
        \label{eq:concgenf}
        PC_r(t) =  \sum_{n=0}^\infty p_{c}(n,r)t^n =
       \frac{Q_r(t)}{\prod_{i=1}^r \left( 1 -
           t^{\frac{i(i+1)}{2}}\right)}
      \end{equation}
      where \(Q_r(1)=1\), and the numerator has degree strictly smaller  than
      \(r^3/6 + r^2/2 + r/3\).
    \item       The proportion of concave partitions with
      at most \(r\) parts among all partitions with at most \(r\) parts
      is
    \begin{equation}
      \label{eq:prop2}
      \frac{r!}{\prod_{i=1}^r \frac{i(i+1)}{2}}.
    \end{equation}
      \end{enumerate}
    \end{cor}
    \begin{proof}
      The only thing which does not follow immediately from
      substituting \(x_i=t\) in the previous theorem is the assertion
      about the degree of the numerator.   From Stanley's ``grey book''
      \cite[Theorem 4.6.25]{Stanley:En1} we have that the rational
      function \(PC_r(t,\dots,t)\) is of degree \(<0\). The degree of
      the denominator is 
      \begin{displaymath}
        \sum_{i=1}^r \frac{i(i+1)}{2} = 
        \frac{r^3}{6} + \frac{r^2}{2} + \frac{r}{3}
      \end{displaymath}
      so the result follows.
    \end{proof}

    We can therefore say with absolute
    certainty  that the first \(Q_r(t)\) are as follows: 
    \begin{equation}
      \label{eq:Qrt}
      \begin{split}
        Q_1(t)  &= 1\\
        Q_2(t)  & = 1+t^2-t^3\\
        Q_3(t) &= 1+t^2+t^5-2t^6-t^8+t^9\\
        Q_4(t) &=
        1 + t^2 + t^4 + t^5-t^6-t^7 + 2t^9-2t^{10}-t^{11}-2t^{12} +
        \\ &+
        2t^{13}-t^{14}-t^{15} + t^{16} + t^{17} + t^{18}-t^{19}
      \end{split}
    \end{equation}
    Hence, we belive that
    \begin{equation}
      \label{eq:PC3}
      PC(t) = \frac{1 + t^2 + O(t^3)}{\prod_{i=1}^\infty \left( 1 -
           t^{\frac{i(i+1)}{2}}\right)}
    \end{equation}

    We've calculated that
    \begin{multline}
      \label{eq:PC2}
      PC(t)=\sum_{n=0}^\infty p_c(n)t^n =
      1+ t + 2 t^2 + 3 t^3 + 4 t^4 + 7 t^5 + 9 t^6 + 11 t^7 +
      \\ +17t^8  +23 t^9 + 28 t^{10} + 39 t^{11} + 48 t^{12}
      + 59 t^{13} + 79 t^{14} +
      \\ +      100 t^{15} + 121 t^{16}
      +152 t^{17} + 185 t^{18} + 225 t^{19} + 280 t^{20} + O(t^{21})
    \end{multline}
    It seems likely that \(\log p_c(n)\) grows as
    \(n^{1/3}\) (i.e. approximately as fast as pseudo-convex
    partitions), but we can not prove this, since we have no estimates
    of the numerator in \eqref{eq:PC3}.

    We have submitted \(\bigl(p_c(n)\bigr)_{n=0}^\infty\) to the OEIS \cite{Sloane};
    it is 
    \htmladdnormallink{A084913}{http://www.research.att.com/projects/OEIS?Anum=084913}.
    The sequences \(\bigl(p_c(n,r)\bigr)_{n=0}^\infty\)
    are 
    \htmladdnormallink{A086161}{http://www.research.att.com/projects/OEIS?Anum=086161},   
    \htmladdnormallink{A086162}{http://www.research.att.com/projects/OEIS?Anum=086162},   
    and 
    \htmladdnormallink{A086163}{http://www.research.att.com/projects/OEIS?Anum=086163}   
    for \(r=2,3,4\).

\end{section}
\raggedright
\bibliographystyle{plain}
\bibliography{journals,articles,snellman,newconcbib}
\end{document}